\documentclass[12pt]{amsart}
\usepackage{graphicx}
\usepackage{graphics}
\usepackage{amsmath}
\usepackage{amscd}
\usepackage{latexsym}
\usepackage{subfigure}
\usepackage{caption}
\usepackage{hyperref}

\begin{document}

\textwidth 5.9in
\textheight 7.9in

\evensidemargin .75in
\oddsidemargin .75in

\newtheorem{Thm}{Theorem}
\newtheorem{Lem}[Thm]{Lemma}
\newtheorem{Cor}[Thm]{Corollary}
\newtheorem{Prop}[Thm]{Proposition}
\newtheorem{Rm}{Remark}

\def\a{{\mathbb a}}
\def\C{{\mathbb C}}
\def\A{{\mathbb A}}
\def\B{{\mathbb B}}
\def\D{{\mathbb D}}
\def\E{{\mathbb E}}
\def\R{{\mathbb R}}
\def\P{{\mathbb P}}
\def\S{{\mathbb S}}
\def\Z{{\mathbb Z}}
\def\O{{\mathbb O}}
\def\H{{\mathbb H}}
\def\V{{\mathbb V}}
\def\Q{{\mathbb Q}}
\def\Cn{${\mathcal C}_n$}
\def\CM{\mathcal M}
\def\CG{\mathcal G}
\def\CH{\mathcal H}
\def\CT{\mathcal T}
\def\CF{\mathcal F}
\def\CA{\mathcal A}
\def\CB{\mathcal B}
\def\CD{\mathcal D}
\def\CP{\mathcal P}
\def\CS{\mathcal S}
\def\CZ{\mathcal Z}
\def\CE{\mathcal E}
\def\CL{\mathcal L}
\def\CV{\mathcal V}
\def\CW{\mathcal W}
\def\IC{\mathbb C}
\def\IF{\mathbb F}
\def\IK{\mathcal K}
\def\IL{\mathcal L}
\def\IP{\bf P}
\def\IR{\mathbb R}
\def\IZ{\mathbb Z}

\title{Cork twisting Schoenflies problem }
\author{Selman Akbulut}
\thanks{Partially supported by NSF grants DMS 0905917}
\address{Department  of Mathematics, Michigan State University,  MI, 48824}
\email{akbulut@math.msu.edu }
\date{\today}
\begin{abstract} 
The stable Andrews-Curtis conjecture in combinatorial group theory is  the statement that every balanced presentation of the trivial group  can be simplified to the trivial form by elementary moves corresponding to ``handle-slides" together with ``stabilization" moves. Schoenflies conjecture is the statement that the complement of any smooth embedding $S^{3}\hookrightarrow S^{4}$  are pair of smooth balls. Here we suggest an approach to these problems by certain cork twisting operation on contractible manifolds, and demonstrate it on the example of the first  Cappell-Shaneson homotopy sphere.
\end{abstract}

\date{}
\maketitle

\setcounter{section}{-1}

\vspace{-.3in}

\section{Introduction}

Let $G(P)=\{ x_1,x_2,..,x_n \;|\; r_{1}(x_1,..,x_n),..,r_{n}(x_1,..x_n) \}$ be a balanced presentation $P$ of the trivial group. Here balanced means  the presentation has the same number of generators and relators. When there is no danger of confusion we will  abbreviate $r_{j}:=r_{j}(x_1,..,x_n)$. The presentation $P$ is called  {\it stably Andrews-Curtis trivial} (SAC-trivial in short) if by changing  relators by the following finite number of the steps, and their inverses, we obtain the trivial presentation: 
\begin{itemize}
\item[(a)] $r_{i}\mapsto r_{i}r_{j} $ for some $j\neq i $.
\item[(b)] Add a new generator $x_{n+1}$ and a relation  $x_{n+1}\gamma$.
\item[(c)] $r_{i} \mapsto r_{i}^{-1}$ or 
$r_{i} \mapsto \gamma r_{i}\gamma^{-1}$, where $\gamma$ represents any word in $G(P)$.
\end{itemize}

\vspace{.05in}

Fundamental group of any compact $2$-complex gives such a presentation, so any compact contractible $4$-manifold $W$, which is a {\it $2$-handlebody} (i.e.\;a handlebody  consisting of handles of index $\leq 2$) has such a presentation. Generators $\{x_{j} \}$ correspond to the $1$-handles, and the relations $\{r_{j}\}$ correspond to the $2$-handles. (a) corresponds to sliding $2$-handles over each other, and (b) corresponds to introducing (or taking away) a canceling pair of $1$ and $2$-handles. Call a pair of $2$-handlebodies SAC {\it equivalent} if they are related by these two steps. 

\vspace{.05in}

Not much known about which presentations of the trivial group are SAC-trivial. in \cite{ak1} the following examples were proposed (n=0,1,..)

$$G(P_{n})=\{ x,y \;|\; xyx=yxy, \; x^{n+1}= y^{n}  \}$$ 

$G(P_{n})$ is the trivial group since the first relation gives $y = (yx)^{-1}x (yx)$, so $y^{n+1}=(yx)^{-1}x^{n+1} (yx)= (yx)^{-1}y^{n} (yx)=x^{-1}y^{n}x=x^{-1}x^{n+1}x=y^{n}$, hence $y=1$ and $x=1$. Gersten showed that $G(P_{2})$ is SAC-trivial (\cite{ge}, \cite{gst}), but it is not known whether the other $G(P_{n})$ are SAC-trivial. The case $G(P_{4})$ is particularly interesting, since it is the fundamental group  presentation of the $2$- handlebody of a Cappell-Shaneson homotopy $4$-ball $W_{0}= \Sigma_{0}-B^{4}$ constructed in \cite{ak1}, where $\Sigma_{0}$ is the $2$-fold covering space of the first Cappell-Shaneson exotic $\R\P^4$  defined in \cite{cs}.

\vspace{.05in}

\begin{figure}[ht]  \begin{center}  
\includegraphics[width=.5\textwidth]{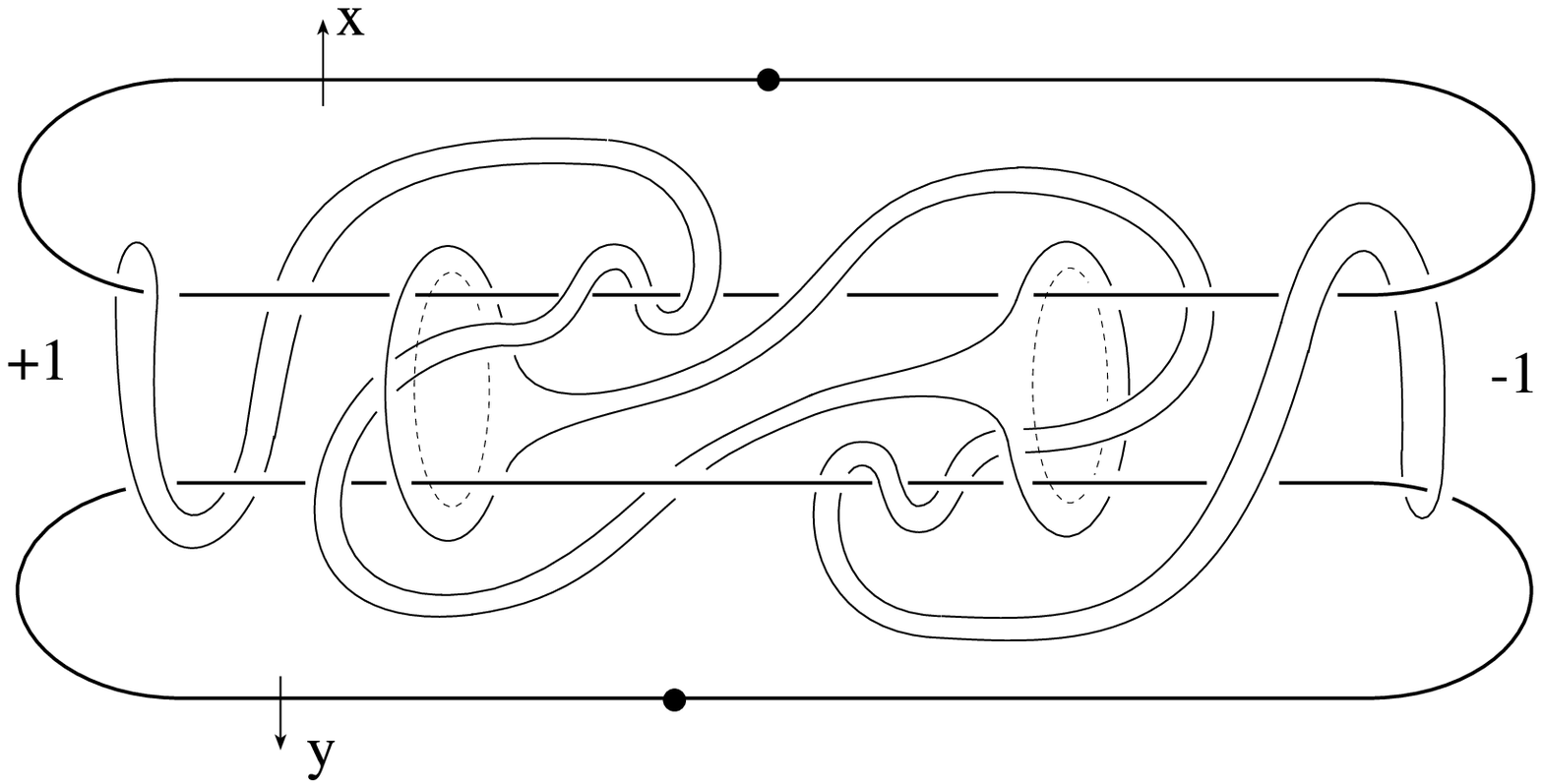}   
\caption{$W_{0}$} 
\label{fig1}
\end{center}
\end{figure}

The main reason topologists are interested SAC problem is its relation to the Schoenflies  conjecture, which says ``{\it The complement of a smoothly imbedded $S^{3}\hookrightarrow S^{4}$ is a disjoint union of two smooth $4$-balls}''. So far, in dimension $4$ only the topological version of this conjecture is known (\cite{b}, \cite{m}). If  the presentation of $\pi_{1}(W)$ a smooth contractible  $2$-handlebody $W^{4}$ is SAC-trivial then $W\times [0,1]=B^{5}$ (because in dimension $5$ canceling handles algebraically is equivalent canceling them geometrically), hence this gives an imbedding $W\hookrightarrow S^{4}$ via its double 
$$D(W):=W\cup_{\partial} -W=\partial (W\times [0,1])=S^4$$ Then since $\partial W=S^{3}$ the topological Schoenflies theorem implies $W$ is homeomorphic to $B^{4}$. To apply this  to the smooth Poincare conjecture (PC), we first puncture a given smooth homotopy $4$-sphere $\Sigma$ to a homotopy $4$-ball $W=\Sigma -B^{4}$ and turn it to $2$-handlebody by canceling its $3$-handles (if we can), reducing it to a SAC problem. This is what is done for the first Cappell-Shaneson homotopy $4$-ball $W_{0}$ in \cite{ak1} (similar proof for the other ones), but there the associated SAC problem was bypassed by directly imbedding $W_{0}\hookrightarrow S^{4}$, hence reducing the smooth PC to Schoenflies problem. In particular if the complement of $W$ in $S^4$ is $B^4$ then $W$ itself must be $B^{4}$.
Figure~\ref{fig1} shows the $3$-handle free handlebody picture of $W_{0}$ which imbeds into $S^{4}$ (\cite{ak1} Figure 28).

In the end it turned out that all Cappell-Shaneson homotopy balls $W_{0}, W_{1},..$ are diffeomorphic to $B^{4}$; without even appealing to the Schoenflies problem (\cite{g}, \cite{a1}). More specifically, the proofs proceed by first introducing canceling 2/3 handle pairs, then canceling all the handles ending up with $W_{i}=B^4$, $i=0,1,2..$. 
Here we revisit the Schoenflies problem by analyzing  the approach of \cite{ak1} more closely, where another a  $2$-handlebody $W_{0}^{*}$ with $\partial W^{*}_{0}=S^3$ was constructed so that $$S^{4}=W_{0}\cup_{\partial} -W_{0}^{*}$$
We can reduce this Schoenflies problem to another Schoenflies problem which we can solve, i.e. by imbedding $W_{0}^{*}\hookrightarrow S^{4}$ with complement $B^{4}$
  $$S^{4}= W_{0}^{*}\cup_{\partial} B^4$$ 
which implies  $W_{0}^{*}=B^{4}$, and so $W_{0}=B^4$. Of course this last step is not new, it is just a case of proving some homotopy $4$-balls are standard by introducing a single canceling $2/3$ handle pair (\cite{g}, \cite{a1}). We  stated it this way to relate it to Schoenflies problem. Curiously the associated presentation of the fundamental group of $W_{0}^{*}$ is $G(P_{2})$, while $\pi_{1}(W_{0})$ is $G(P_{4})$. The hope is, associating $W$ another convenient ``{\it twin}" contractible manifold $W\mapsto W^{*}$ might help to resolve SAC triviality. 
   
\section{Flexible contractible manifolds} A flexible contractible $4$-manifold is a smooth compact contractible $2$-handlebody, where its $2$-handles are represented by $0$-framed unknotted curves (i.e. after erasing circles with dots we get a $0$-framed unlink).

 \begin{figure}[ht]  \begin{center}  
\includegraphics[width=.53\textwidth]{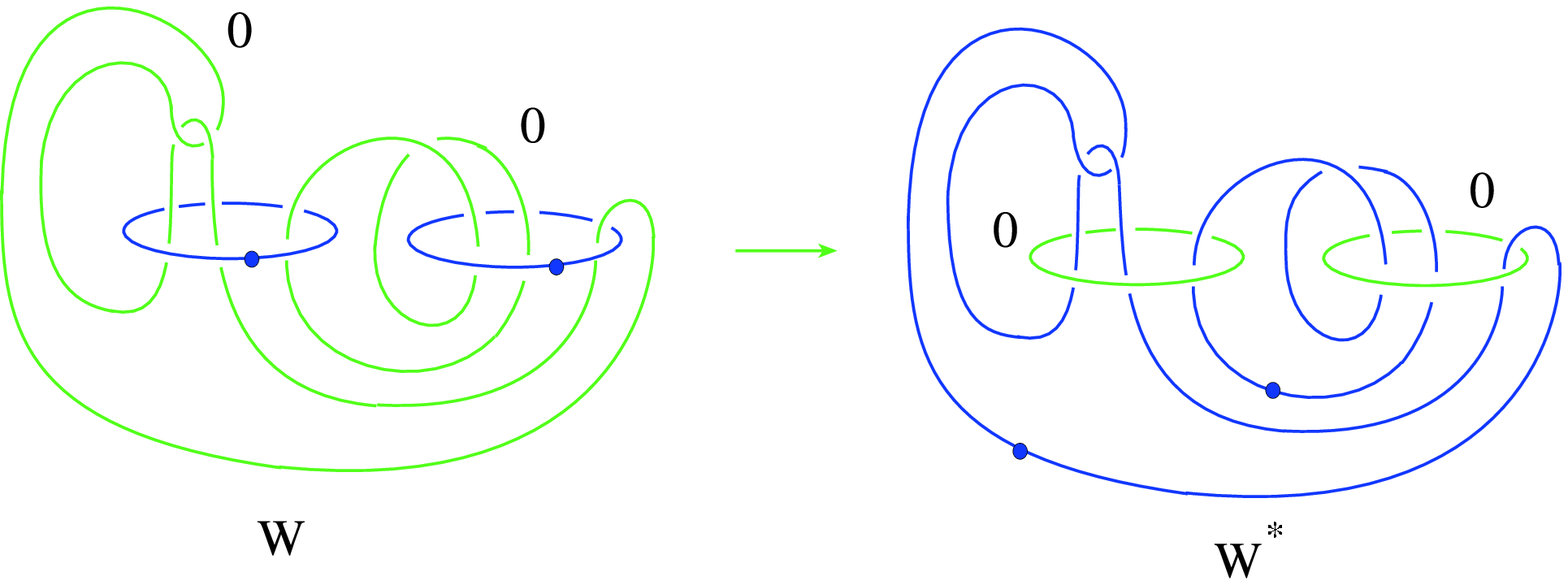}   
\caption{A flexible contactible manifold and its twin} 
\label{fig2}
\end{center}
\end{figure} 

We call the $2$-handlebody $W^{*}$ obtained from $W$ by zero and dot exchanges of its handles ($S^{2}\times B^2 \leftrightarrow B^{3}\times S^{1}$ exchanges in  the interior of $W$) the {\it twin of W} . We call the operation $W\mapsto W^{*}$ {\it cork twisting the flexible manifold $W$}. Notice that this notion depends on the handles. Here we do not address the problem of how the twin of $W$ changes after handle slides of W. It is clear that  this  decomposes the $4$-sphere as $S^{4}= W\cup_{\partial} -W^{*}$, i.e. $W$ imbeds into $S^4$ with  complement $W^{*}$.

\begin{figure}[ht]  \begin{center}  
\includegraphics[width=.34\textwidth]{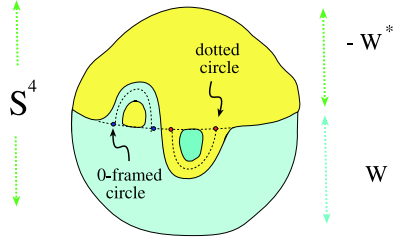}   
\caption{A flexible contractible manifold and its twin} 
\label{fig3}
\end{center}
\end{figure}

\section{The twin of $W_{0}$}

\begin{figure}[ht]  \begin{center}  
\includegraphics[width=.5\textwidth]{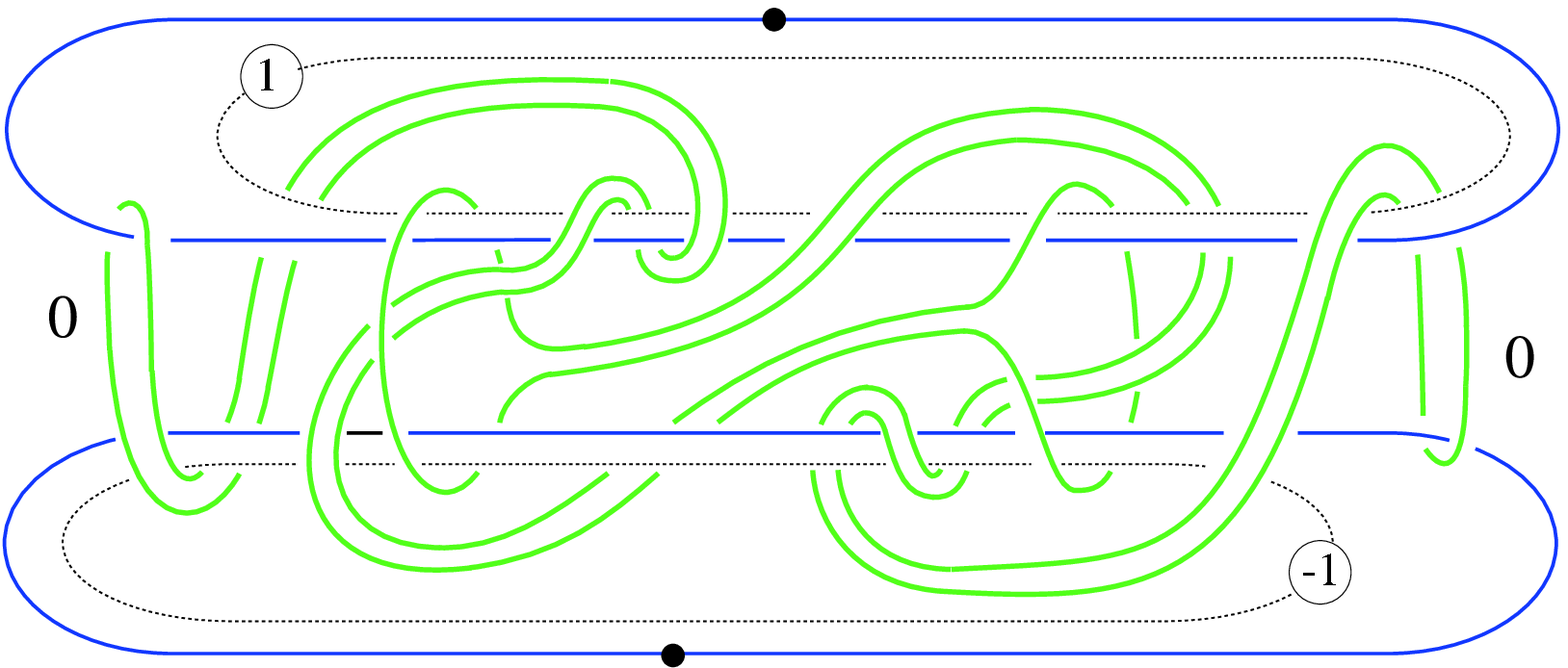}
\caption{$W_{0}$}    
\label{fig4}
\end{center}
\end{figure}

First  by applying the diffeomorphism 
described in Section 2.3 of \cite{a3}, we identify the handlebody $W_{0}$ of Figure~\ref{fig1} with  Figure~\ref{fig4}. This diffeomorphism is a combination of introducing and canceling 1/2 handle pairs. Recall that the associated presentation of $\pi_{1}(W_{0})$ is $G(P_{4})$.

\begin{figure}[ht]  \begin{center}  
\includegraphics[width=.35 \textwidth]{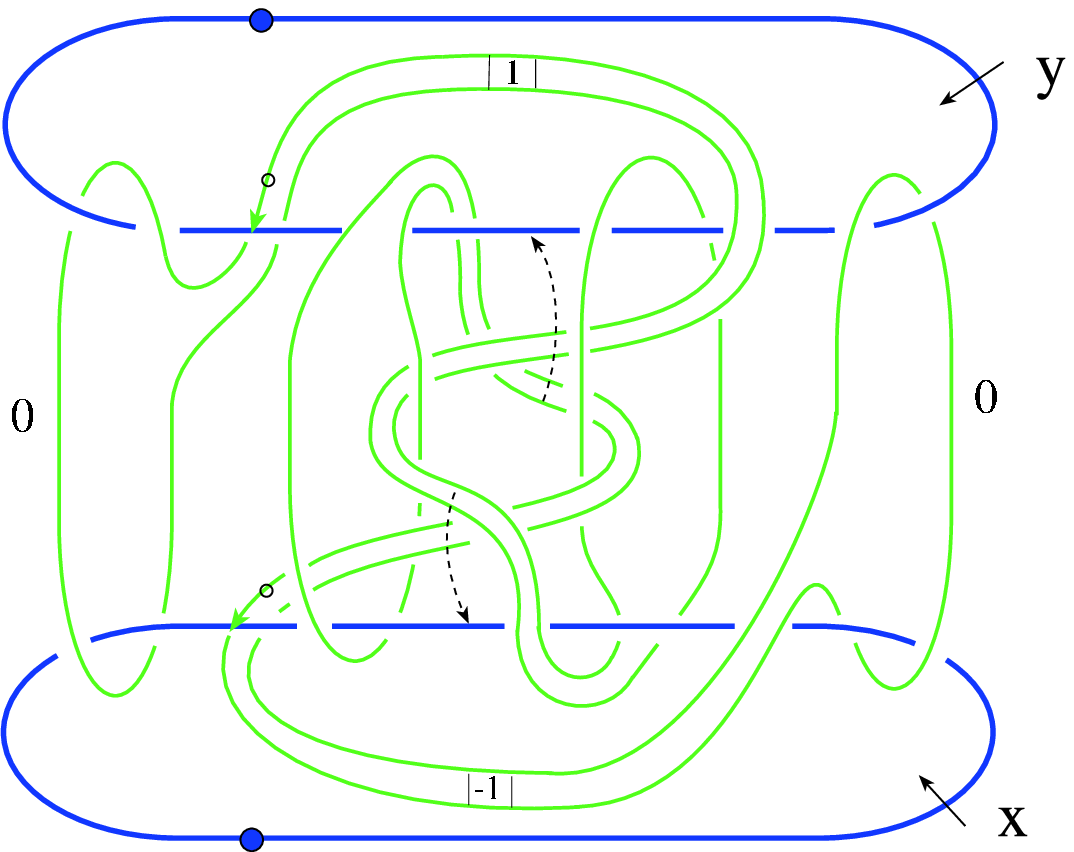}   
\caption{$W_{0}^{*}$} 
\label{fig5}
\end{center}
\end{figure}

\begin{Prop}
The twin $W_{0}^{*}$ of the $2$-handlebody of $W_{0}$ in Figure~\ref{fig4} is  given by Figure~\ref{fig5}, and the associated presentation of $\pi_{1}(W_{0}^{*})$ is  $G(P_{2})$.
\end{Prop}

\begin{figure}[ht]  \begin{center}  
\includegraphics[width=.5 \textwidth]{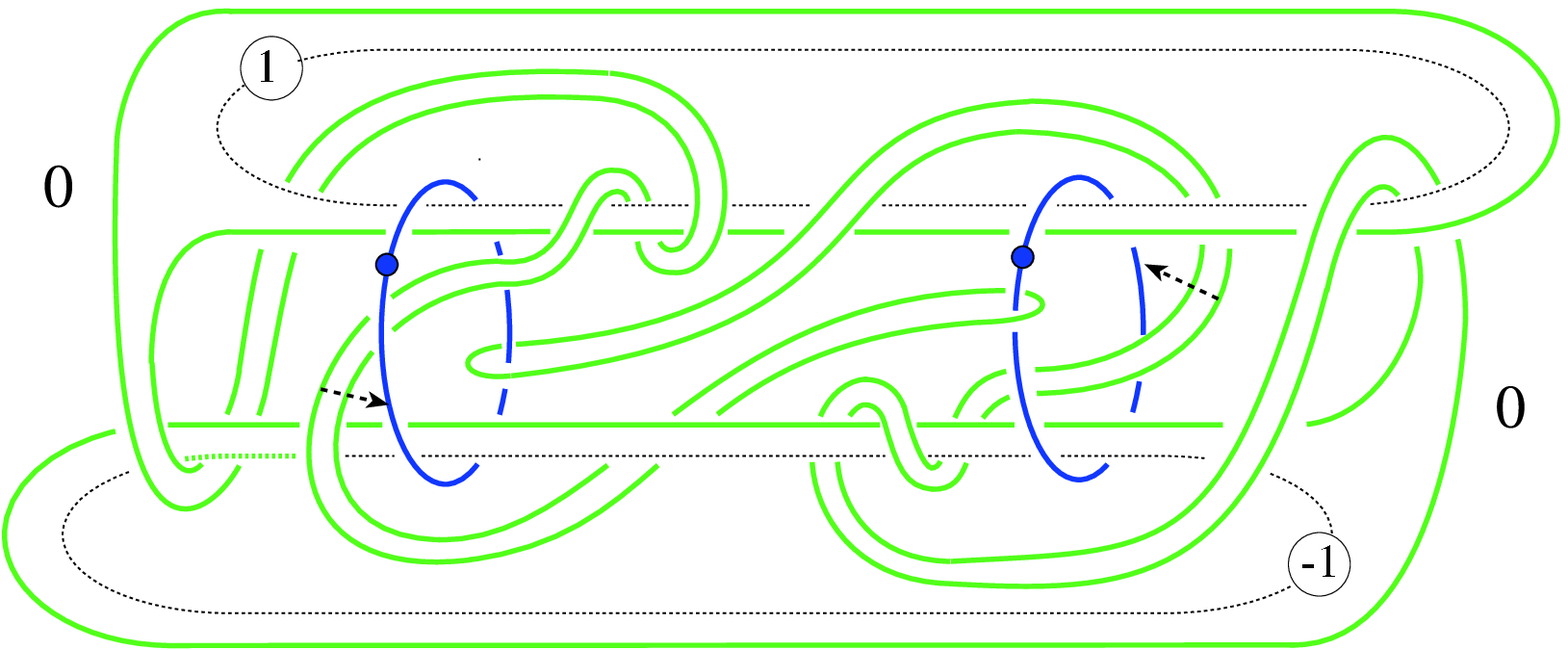}   
\caption{$W_{0}^{*}$} 
\label{fig6}
\end{center}
\end{figure}

\proof Clearly Figure~\ref{fig6} is the twin of $W_{0}$ in Figure~\ref{fig4}. Notice that the two circles with dots do not link each other, since the affect of blowing down one $+1$ is undone by blowing down another $-1$ (before putting dots on them),  of course during this process the framed circles representing the $2$-handles changed appropriately.

\begin{figure}[ht]  \begin{center}  
\includegraphics[width=.37 \textwidth]{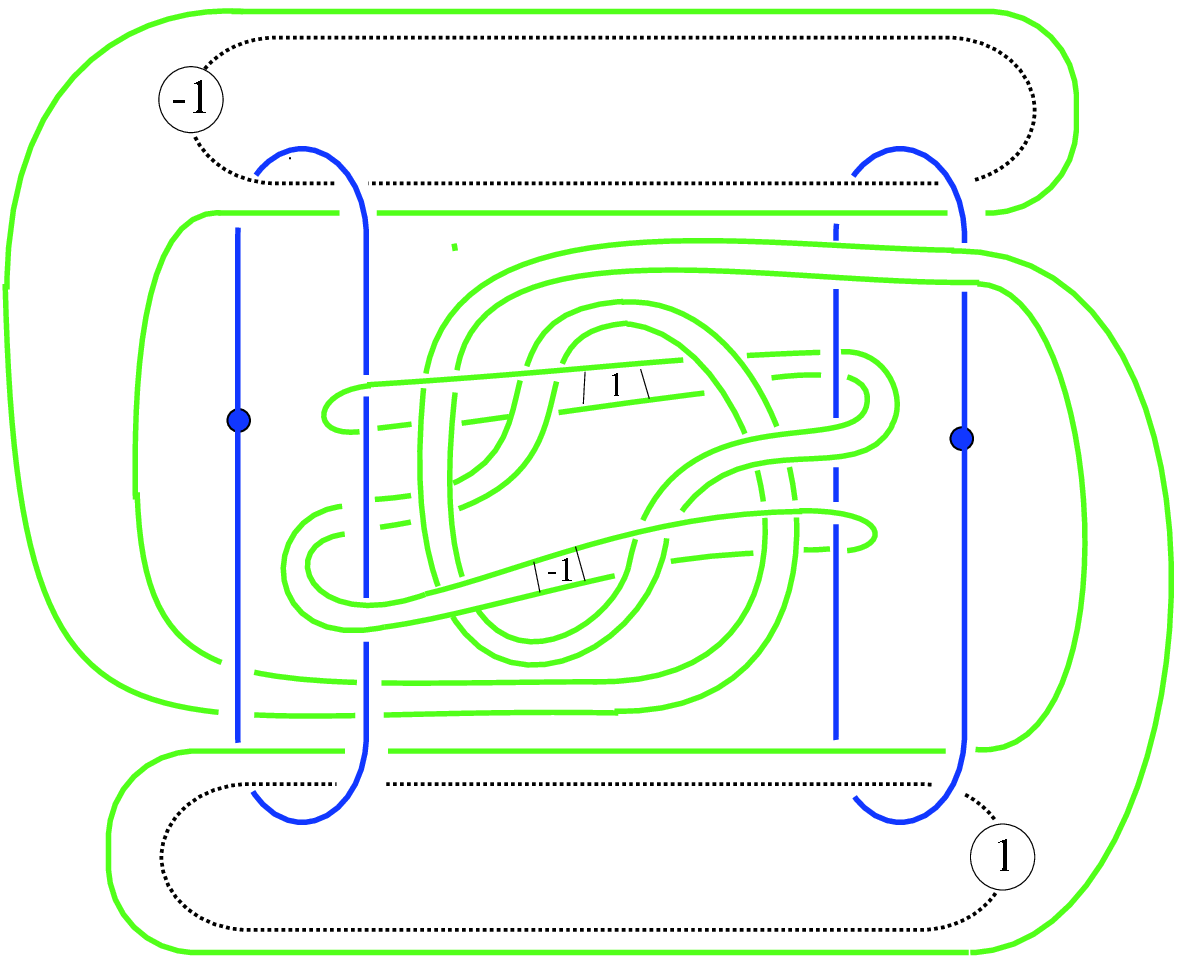}   
\caption{$W_{0}^{*}$} 
\label{fig7}
\end{center}
\end{figure}

By sliding $2$-handles over the $1$-handles (as indicated by dotted arrows) in Figure~\ref{fig6} we get Figure~\ref{fig7}. Then isotoping the two dotted circles away from each other and after rotating the figure $90^{0}$ we get Figure~\ref{fig5}.

Next we calculate the presentation of $\pi_{1}(W_{0}^{*})$  from Figure~\ref{fig5}:
\begin{itemize}
\item[(a)] $x^{2}yx^{-1}yx^{-1}y^{-1}=1$
\item[(b)] $y^{-2}x^{-1}yx^{-1}yx=1$
\end{itemize}
(a) $\implies x^{-1}yx y^{-1}=xyx^{-1}$, and (b) $\implies x^{-1}yx y^{-1}=y^{-1}xy$. Hence $xyx=yxy$. Also (a) 
$\implies  x^{3}yx^{-1}y=xyx=yxy \implies x^{3}y=yx^{2}$. Hence
 $x^{3}=yx^{2}y^{-1}=(yxy^{-1})^{2}=(x^{-1}yx)^{2}=x^{-1}y^{2}x\implies x^{3}=y^{2}$.  Hence we get the presentation $G(P_{2})$. \qed
 
 \begin{Rm} As shown in \cite{ak1}, attaching pair of $2$-handles to $\partial W_{0}$ along the dotted circles of Figure~\ref{fig1} gives $\#_{2}(S^{2}\times B^{2})$, which can be capped by a pair of $3$-handles $\#_{2}(B^{3}\times S^{1})$. Reader can check that turning these handle pairs upside down gives the handlebody  $W_{0}^{*}$.
 \end{Rm}
 
 \begin{Prop} $W_{0}^{*}$ smoothly imbeds into $S^{4}$ with complement $B^{4}$, hence $W_{0}^{*} $ (therefore $W_{0}$) is diffeomorphic to $B^{4}$.   
 \end{Prop}
 
\begin{figure}[ht]
\centering
\begin{minipage}{.35\textwidth}
  \centering
  \includegraphics[width=.8\linewidth]{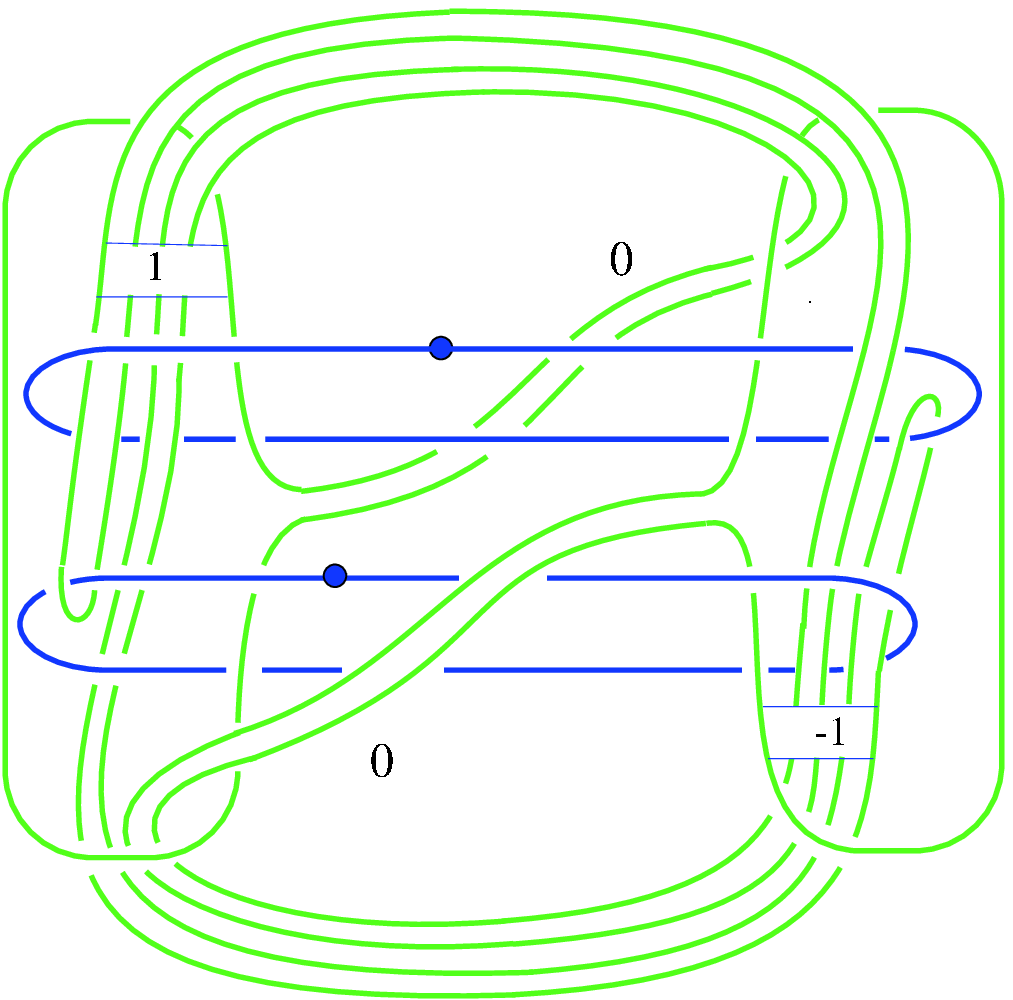}
  \caption{$W_{0}^{*}$}
  \label{fig8}
\end{minipage}%
 \begin{minipage}{.4\textwidth}
  \includegraphics[width=.85\linewidth]{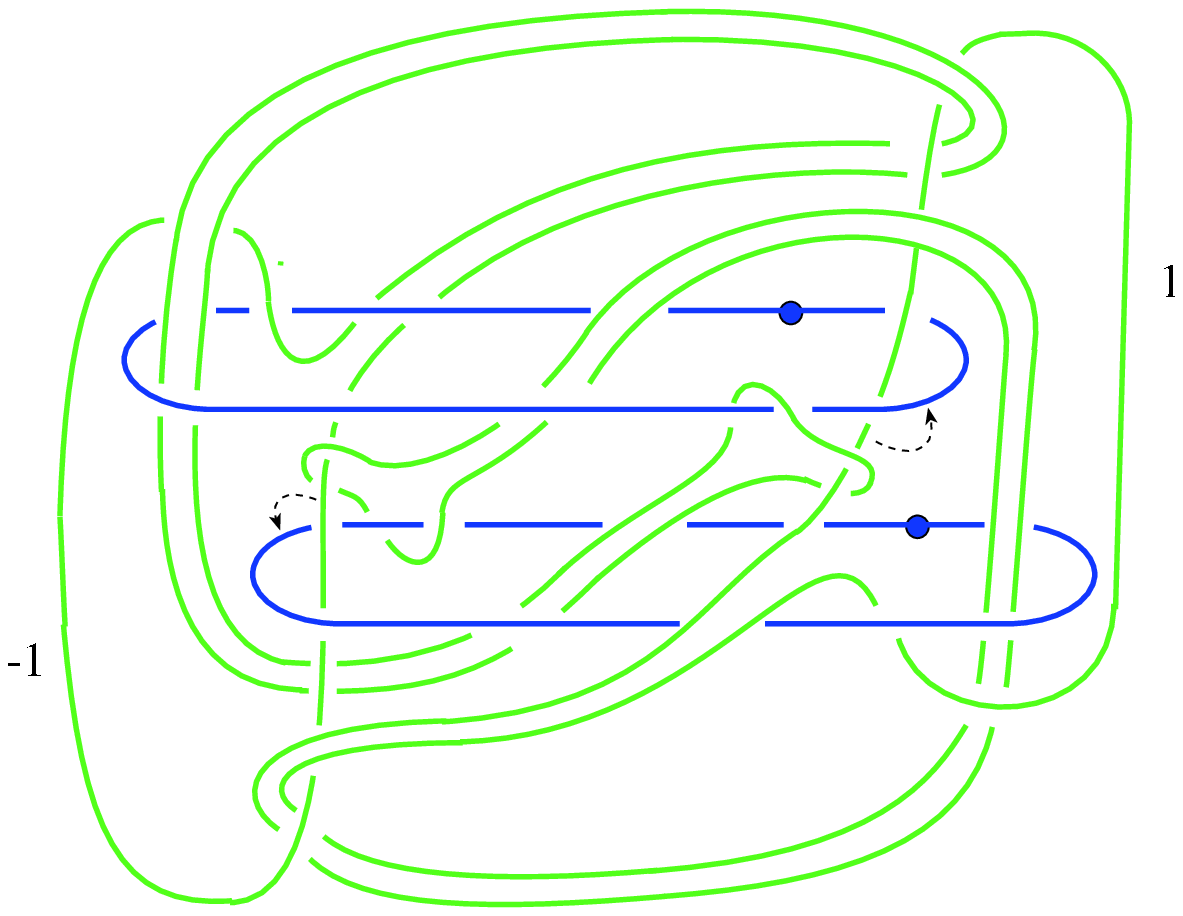}
  \caption{$W_{0}^{*}$}
  \label{fig9}
\end{minipage}
\end{figure}

\begin{figure}[ht]  \begin{center}  
\includegraphics[width=.5\textwidth]{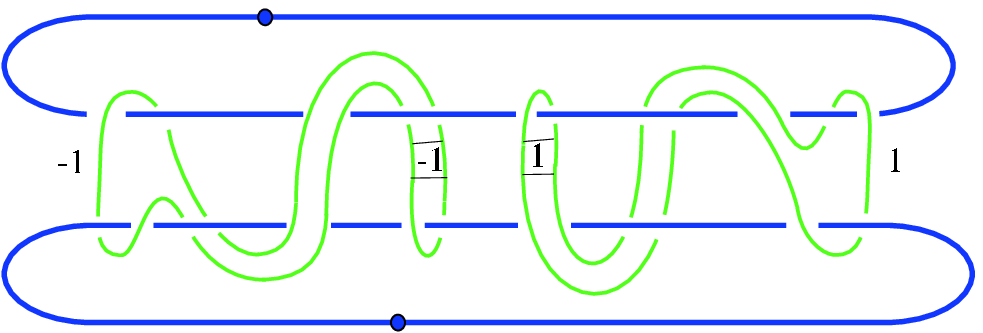}   
\caption{$W_{0}^{*}$} 
\label{fig10}
\end{center}
\end{figure}

\begin{figure}[ht]
\centering
\begin{minipage}{.4\textwidth}
  \centering
  \includegraphics[width=.7\linewidth]{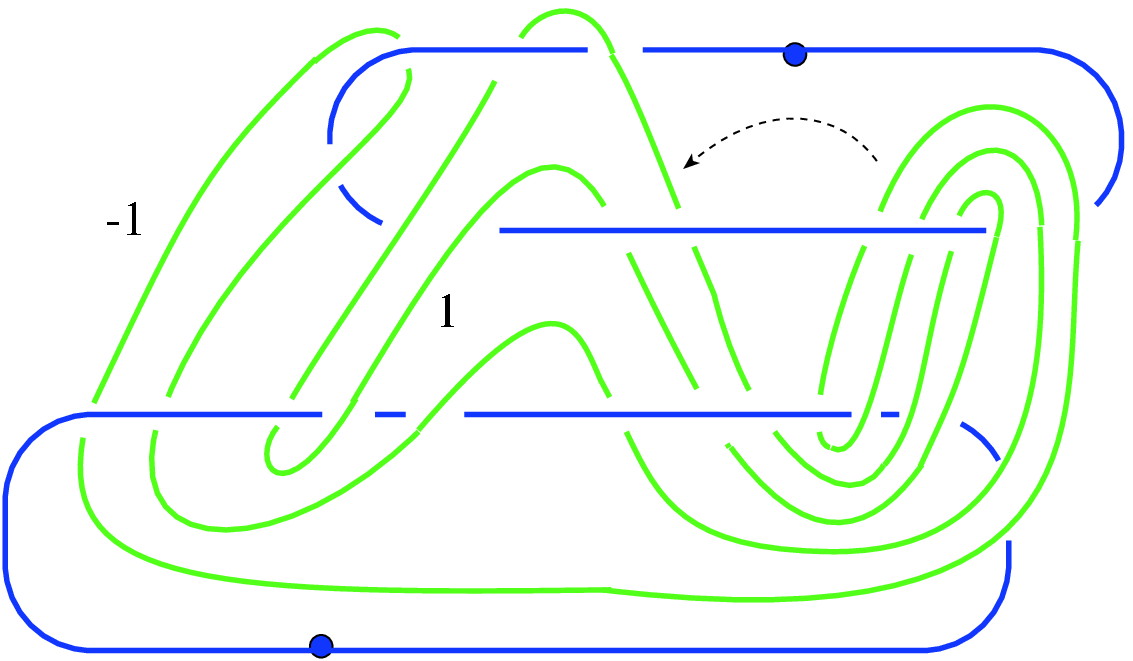}
  \caption{$W_{0}^{*}$}
  \label{fig11}
\end{minipage}%
 \begin{minipage}{.45\textwidth}
  \includegraphics[width=.7\linewidth]{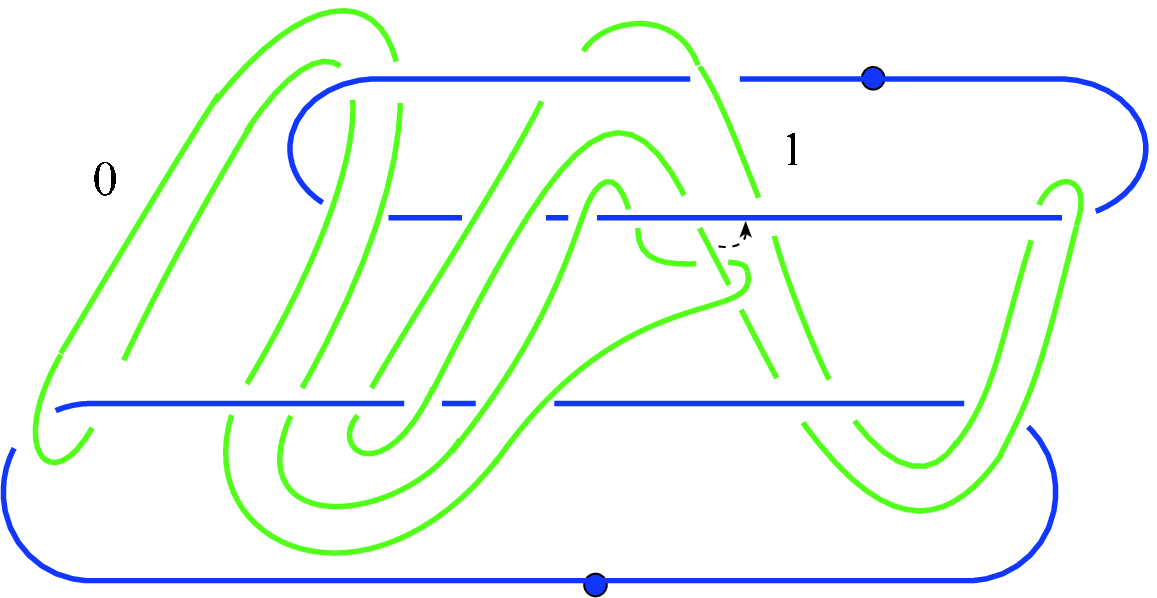}
  \caption{$W_{0}^{*}$}
  \label{fig12}
\end{minipage}
\end{figure}

 \proof By sliding the $2$-handles over the $1$-handles of Figure~\ref{fig5} (as indicated by the dotted arrows in the figure) we get Figure~\ref{fig8}. Again by applying the diffeomorphism of Section 2.3 of [A3] to this figure twice, we get Figure~\ref{fig9}, then after the  handle slides (indicated by dotted arrows) we get  Figure~\ref{fig10}, and then by an isotopy we get Figure~\ref{fig11}.

\begin{figure}[ht]  \begin{center}  
\includegraphics[width=.65 \textwidth]{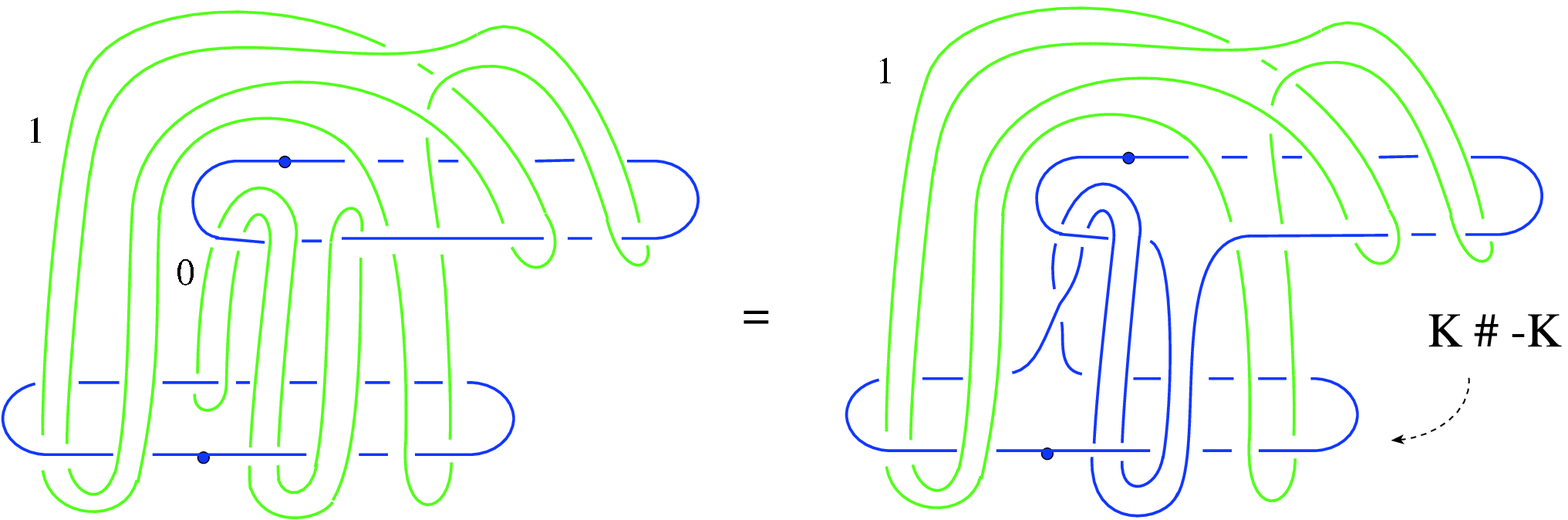}   
\caption{$W_{0}^{*}$} 
\label{fig13}
\end{center}
\end{figure}

 We  will now modify Figure~\ref{fig11} by a sequence of handle slides, and  adding (and canceling) $1/2$ handle pairs to get a new $2$-handlebody presentation of $W_{0}^{*}$ which will have the required property. First by doing the handle slides indicated by the dotted arrows, we go from Figure~\ref{fig11} to Figures~\ref{fig12} and ~\ref{fig13}. Then by canceling a $1/2$ handle pair we obtain the second picture of Figure~\ref{fig13}, which is a ribbon $1$-handle, induced from $K\#-K$ where $K$ is the trefoil knot,  and a $2$-handle.


\begin{figure}[ht]  \begin{center}  
\includegraphics[width=.56 \textwidth]{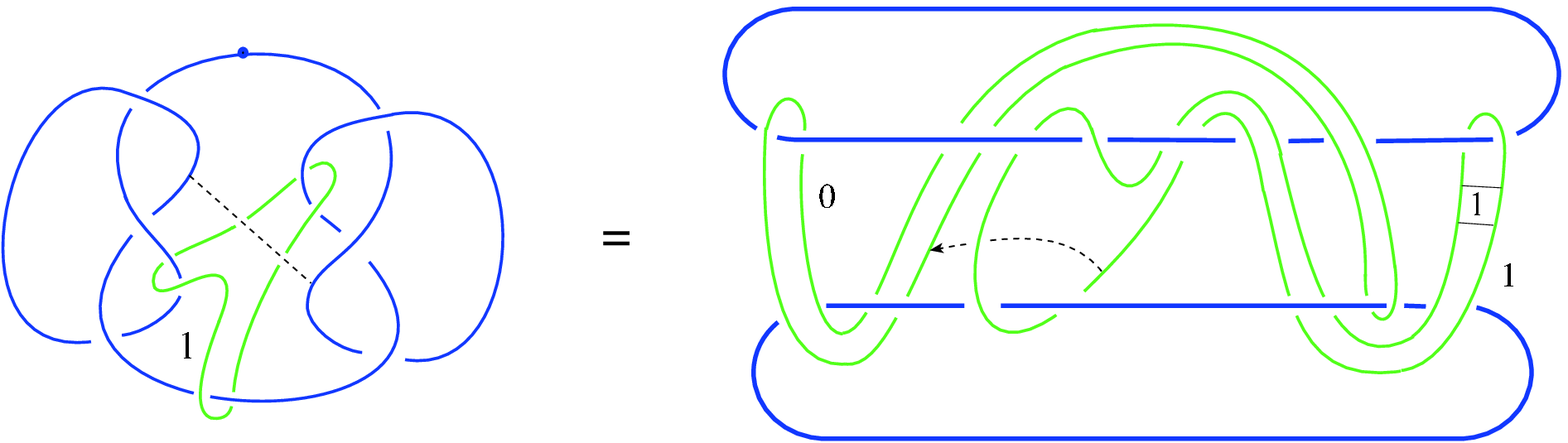}   
\caption{$W_{0}^{*}$} 
\label{fig14}
\end{center}
\end{figure}

\begin{figure}[ht]  \begin{center}  
\includegraphics[width=.7 \textwidth]{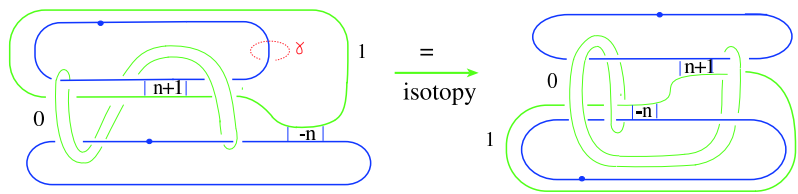}   
\caption{$C_{n}$ and a useful isotopy} 
\label{fig15}
\end{center}
\end{figure}

\noindent The first picture of Figure~\ref{fig14} is the short hand of this handlebody (dotted line indicates the ribbon move giving the ribbon $1$-handle in Figure~\ref{fig13}). The second picture of Figure~\ref{fig14} is drawn after this ribbon move. Then doing the indicated handle slide to Figure~\ref{fig14} we get the handlebody $C_{2}$, where $C_{n}$ is the handlebody described in Figure~\ref{fig15}. 

Notice  $C_{n}$ gives the presentation $G(P_{n})$ and is similar to $H_{n,1}$ of \cite{g}, but it differs in the way its $0$-framed $2$-handle links the $1$-handles (this fact provides us the useful isotopy of Figure~\ref{fig15}).  This difference is due to the fact that here we are getting a nonstandard ribbon which $K\#-K$ bounds.  Now the proof of the proposition follows from the following Lemma~\ref{lem3} whose proof is similar  to the one in \cite{g} for $H(n,1)$. \qed.

\begin{Lem}\label{lem3}
 $C_{n}$ smoothly imbeds into $S^{4}$ with complement $B^{4}$. 
\end{Lem}

\proof Attach a $2$-handle to $C_{n}$ along the loop 
$\gamma $ with $-1$ framing, as shown  in the first picture of Figure~\ref{fig15} (this framing corresponds to $0$-framing when viewed from  $S^{3}$). Denote this manifold by $C_{n} +\gamma^{-1}$. The steps of Figure~\ref{fig15} show the following equivalences by handle slides $$C_{n} +\gamma^{-1} \cong C_{n-1} +\gamma^{-1}...\cong C_{0}+\gamma^{-1}$$

 \begin{figure}[ht]  \begin{center}  
\includegraphics[width=.7 \textwidth]{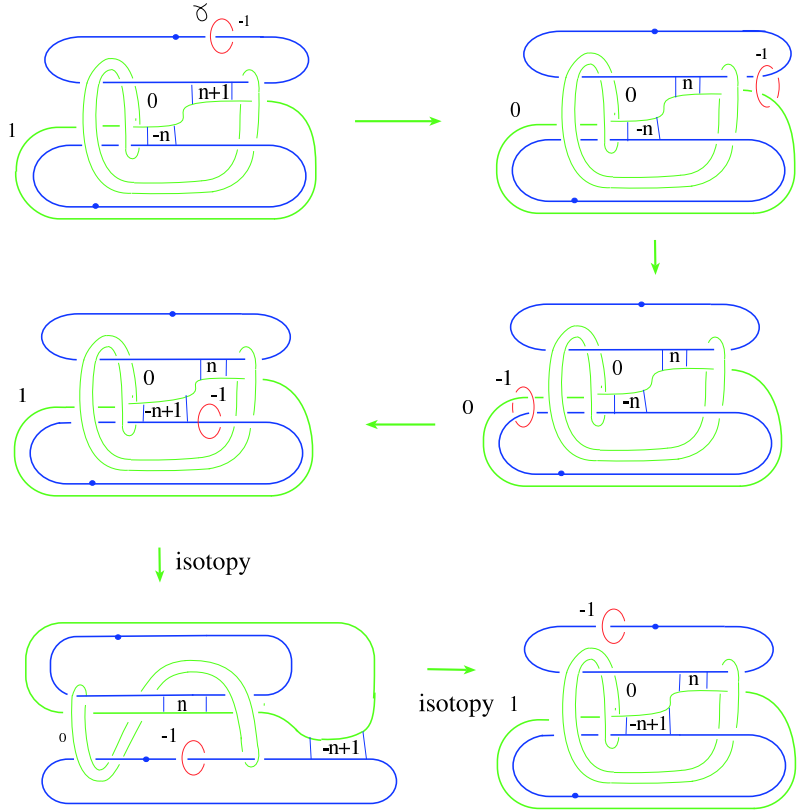}   
\caption{$C_{n}+ \gamma^{-1} \cong C_{n-1}+ \gamma^{-1}$} 
\label{fig16}
\end{center}
\end{figure}

\newpage

\begin{figure}[ht]  \begin{center}  
\includegraphics[width=.57 \textwidth]{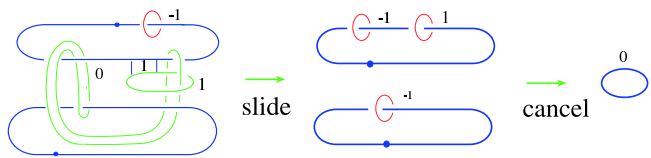}   
\caption{$C_{0}+ \gamma^{-1}\cong S^{2}\times B^{2}$} 
\label{fig17}
\end{center}
\end{figure}

Furthermore, Figure~\ref{fig17} shows $C_{0} + \gamma^{-1} \cong S^{2}\times B^{2}$, hence we can cap $C_{n} + \gamma^{-1}$ with $S^{1}\times B^{3}$ to get $S^{4}$. Let $N$ be the handlebody consisting of $S^{1}\times B^{3}$  union the dual of the $2$-handle $\gamma^{-1}$. $N$ is a contractible manifold with boundary $S^{3}$, consisting of a single pair if  $1$- and $2$-handles, hence $N\cong B^{4}$. This fact follows from  Property-R theorem of \cite{ga}. \qed 

\begin{Rm} There is a certain dictionary  relating AC-triviality of a $2$-handlebody to its twin, which we didn't discussed here, opting directly dealing with Schoenflies problem. This is because the circle with dots can slide over each other just like $2$-handles slide over each other (e.g. Section 1.2 of \cite{a3}). We hope to address this in a future paper. 
\end{Rm}

\section*{Acknowledgement}
I want to thank  Peter Teichner and Rob Kirby for stimulating my interest to the problems discussed here.

\end{document}